\documentclass[11pt]{article}
\usepackage{anysize,ifpdf}\marginsize{3.5cm}{3.5cm}{1.3cm}{2cm}
\ifpdf \pdfpagewidth=210mm\pdfpageheight=297mm \fi

\usepackage{amssymb}
\makeatletter
\def\ps@myheadings{\let\@oddfoot\@empty\let\@evefoot\@empty\def\@oddhead{\hbox{}\hfil\footnotesize\footnotesize\thepage}\def\@evenhead{\footnotesize\footnotesize\thepage\hfil\hbox{}}}
\headsep=\baselineskip
\let\origtitle=\title \renewcommand\title[1]{\origtitle{\normalsize\bfseries\uppercase{#1}}}
\let\origauthor=\author \renewcommand\author[1]{\origauthor{\footnotesize\uppercase{#1}}}
\let\origdate=\date \renewcommand\date[1]{\origdate{\footnotesize #1}}
\let\origthebibliography=\thebibliography \renewcommand\thebibliography[1]{\origthebibliography{#1}\footnotesize\itemsep=0cm\parskip=0cm}
\let\orig@item=\@item \def\@item[#1]{\orig@item[\rm #1]}

\renewcommand\section{\@startsection{section}{1}{\z@}{-1.5\baselineskip plus 0.5\baselineskip minus 0.5\baselineskip}{0.75\baselineskip plus 0.5\baselineskip}{\normalsize\scshape\centering}}
\renewcommand\subsection{\@startsection{subsection}{2}{\z@}{-1\baselineskip}{0.5\baselineskip}{\normalsize\itshape\centering}}
\renewcommand\paragraph{\@startsection{paragraph}{4}{\z@}{1\baselineskip}{-0.5em}{\normalsize\bfseries}}
\renewcommand\@seccntformat[1]{\csname the#1\endcsname.\enspace}
\renewcommand\@begintheorem[2]{\trivlist\item[\hskip\labelsep{\bfseries#1 #2.}]\it}
\renewcommand\@opargbegintheorem[3]{\trivlist\item[\hskip\labelsep{\bfseries#1 #2 {\rm(#3)}.}]\it}
\AtBeginDocument{\renewcommand\@listi{\topsep=0cm \parsep=0cm \itemsep=0cm \partopsep=0.5\baselineskip \labelwidth=\leftmargini \advance\labelwidth by -\labelsep \leftmargin=\leftmargini}
   \let\@listii=\@listi \let\@listiii=\@listi }
\makeatother

\newenvironment{varthm*}[1]{\begin{list}{}{\labelwidth=0cm \leftmargin=0cm \listparindent=\parindent}\item[\hspace{\labelsep}\bfseries #1.]\itshape}{\end{list}}

\newcommand\qedsymbol{\frame{\rule[0pt]{0pt}{8pt}\rule[0pt]{8pt}{0pt}}}
\newcommand\proofname{Proof}
\newenvironment{proof}[1][\proofname]{\trivlist\item[\hskip\labelsep{\it #1.}]}{\hspace*{\fill}\qedsymbol\endtrivlist}

\newcounter{comments}

\renewcommand\ge{\geqslant}
\renewcommand\le{\leqslant}
\renewcommand\epsilon{\varepsilon}
\renewcommand\phi{\varphi}

\newcommand\liste[3]{\mbox{$#1_{#2},\dots,#1_{#3}$}}  
\newcommand\righttext[1]{\qquad\mbox{#1}}
\newcommand\eqnref[1]{(\ref{#1})}
\newcommand\eqdef{=_{\rm def}}
\newcommand\eps{\varepsilon}
\newcommand\Q{\mathbb Q}
\newcommand\R{\mathbb R}

\newcommand\lessdivisor{\preccurlyeq}

\begin{document}

\title{A simple proof for the existence of Zariski decompositions on surfaces}
\author{Thomas Bauer}
\date{November 7, 2007}
\maketitle
\thispagestyle{empty}

   In his fundamental paper \cite{Zar62}, Zariski established the
   following result:

\begin{varthm*}{Theorem}
   Let $D$ be an effective $\Q$-divisor on a smooth projective
   surface~$X$. Then there are uniquely determined effective
   (possibly zero) $\mathbb Q$-divisors $P$ and $N$ with
   $$
      D = P + N
   $$
   such that
\begin{itemize}
\item[(i)]
   $P$ is nef,
\item[(ii)]
   $N$ is zero or has negative definite intersection matrix,
\item[(iii)]
   $P\cdot C=0$ for every irreducible component $C$ of $N$.
\end{itemize}
\end{varthm*}

   The decomposition $D=P+N$
   is called the \textit{Zariski decomposition} of $D$,
   the divisors $P$ and $N$ are respectively the
   \textit{positive} and \textit{negative} parts of $D$.
   Zariski's result
   has been used to study linear series on surfaces, and in the
   classification of surfaces (see \cite[Chapt.~14]{Bad99} and
   \cite[Sect.~2.3.E]{PAG}, as well as the references therein).
   We also mention that
   there is an extension
   to pseudo-effective divisors due to Fujita
   (see \cite{Fuj79} and the nice account in \cite{Bad99}).

   Given an effective divisor $D$,
   Zariski's original proof employs a
   rather sophisticated procedure to construct
   the negative part $N$
   out of those
   components $C$ of $D$
   satisfying $D\cdot C\le 0$.
   Our purpose here is to provide a quick and simple
   proof,
   based
   on the idea that the positive part $P$
   can be constructed
   as a maximal nef subdivisor of $D$.
   This maximality condition is
   in the surface case equivalent to
   the defining condition of Nakayama's
   $\nu$-decomposition of pseudo-effective $\R$-divisors
   (see the Remark below).
   It may be useful that this approach
   yields a practical algorithm for the computation of
   $P$.

\paragraph*{\it Notation.}
   For $\Q$-divisors $P$ and $Q$ we will write $P\lessdivisor Q$, if
   $P$ is a subdivisor of $Q$, i.e., if
   the difference
   $Q-P$ is effective or zero.
   Similarly, we will use the partial ordering $\lessdivisor$
   in $\Q^r$ that is defined by
   $(x_1,\dots,x_r)\lessdivisor(y_1,\dots,y_r)$, if $x_i\le y_i$
   for all $i$.

\begin{proof}[Proof of existence]
   Write $D=\sum_{i=1}^r a_iC_i$ with distinct irreducible curves
   $C_i$ and positive rational numbers $a_i$.
   Consider now all effective $\Q$-subdivisors $P$ of
   $D$, i.e., all divisors of the form
   $P=\sum_{i=1}^r x_iC_i$ with rational coefficients $x_i$
   satisfying
   $0\le x_i\le a_i$.
   A divisor $P$ of this kind is
   nef if and
   only if
   \begin{equation}\label{nefness-condition}
      \sum_{i=1}^r x_i\,C_i\cdot C_j\ge 0 \quad\mbox{ for $j=1,\dots r$.}
   \end{equation}
   This system of linear inequalities for the rational numbers
   $x_i$ has a maximal solution (with respect to
   $\lessdivisor$)
   in the rational
   cuboid
   $$
      [0,a_1]\times\dots\times[0,a_r] \subset\Q^r \ .
   $$
   To see this, note first that the subset $K$ of the cuboid that
   is described by \eqnref{nefness-condition} is a rational
   convex polytope defined by finitely many rational halfspaces.
   It is therefore
   the convex envelope of finitely many rational points.
   We are
   done if $(a_1,\dots,a_r)\in K$. In the alternative case
   consider for rational $t<1$ the family of hyperplanes
   $H_t=\{(x_1,\dots,x_r)\in\Q^r\mid\sum_i x_i=t\sum_i a_i\}$.
   There is then a maximal $t$ such that $H_t$ intersects $K$,
   the point of intersection being a vertex of $K$.

   Let now
   $P=\sum_{i=1}^r b_iC_i$
   be a divisor that is determined by a maximal solution,
   and put $N=D-P$. Then both $P$ and $N$ are effective, and $P$
   is a maximal nef $\Q$-subdivisor of $D$. We
   will now show that (ii) and (iii) are satisfied as well.

   As for (iii): Suppose $P\cdot C>0$ for some component $C$ of
   $N$. As $C\lessdivisor N$, we have $b_i<a_i$, so that for
   sufficiently small rational numbers $\eps>0$, the divisor
   $P+\eps C$ is a subdivisor of~$D$. For curves $C'$ different
   from $C$ we clearly have $(P+\eps C)\cdot C'\ge 0$. Moreover,
   $(P+\eps C)\cdot C=P\cdot C+\eps C^2>0$ for small $\eps$. So
   $P+\eps C$ is nef, contradicting the maximality of $P$.

   As for (ii):
   Supposing that the divisor $N$ is non-zero, we need to show
   that its intersection matrix is negative definite.
   We will prove:

   \begin{itemize}
   \item[(*)]
      If $N$ is a divisor, whose intersection matrix $S$ is not
      negative definite, then there
      is an effective non-zero nef divisor $E$, whose
      components are among those of~$N$.
   \end{itemize}

   Granting (*) for a moment, let us show how to complete the
   proof. Assume by way of contradiction that the intersection
   matrix of $N$ is not negative definite, and take $E$ as in
   (*).
   Consider then for
   rational $\eps>0$
   the $\Q$-divisor
   $$
      P' \eqdef P+\eps E \ .
   $$
   Certainly $P'$ is effective and nef.
   As
   all components of $E$
   are among the components of $N$,
   it is clear that $P'$ is
   a subdivisor of $D$ when $\eps$ is small enough.
   But this is a contradiction, because $P'$
   is strictly
   bigger than
   $P$.

   It remains to prove (*). To this end
   we distinguish between two cases:

\medskip
   \textit{Case 1: $S$ is not negative
   semi-definite.}
   In this case there is a divisor $B$ whose components are
   among those of $N$ such that $B^2>0$.
   Then, writing $B=B'-B''$ as a difference of
   effective divisors
   having no common components, we have
   $
      0<B^2=(B'-B'')^2=B'^2-2B'B''+B''^2,
   $
   and hence $B'^2>0$ or $B''^2>0$. Therefore, replacing $B$ by $B'$
   or $B''$ respectively, we may
   assume that $B$ is effective. But then it follows from
   the Riemann-Roch theorem that the linear series $|\ell B|$ is
   large for $\ell\gg 0$.
   So we can write $\ell B=E_\ell+F_\ell$, where $|E_\ell|$ is the non-zero
   moving
   part of $|\ell B|$. Then $E=E_\ell$ is a nef divisor as
   required, so that the proof of (*) is complete in this case.

\medskip
   \textit{Case 2: $S$ is negative
   semi-definite.} Let $\liste C1k$ be the components of $N$.
   We argue by induction on $k$.
   If $k=1$, then $N^2=C_1^2=0$, so $C_1$ is nef and we are done
   taking
   $E=C_1$. Suppose then
   $k>1$. The hypotheses on $S$ imply that
   $S$ does
   not have full rank. Therefore there is a non-zero divisor $R$,
   whose components are among $\liste C1k$, having the property
   that $R\cdot
   C_i=0$ for $i=1,\dots k$.
   If one of the divisors $R$ or $-R$ is effective, then it is
   nef,
   and we
   are done, taking $E=R$ or $E=-R$ respectively. In the
   alternative case we write $R=R'-R''$, where $R'$ and $R''$ are
   effective non-zero divisors without common components. We have
   $$
      0=R^2=R'^2-2R'R''+R''^2 \ .
   $$
   As by hypothesis $R'^2\le 0$ and $R''^2\le 0$, we must have
   $R'^2=0$. The divisor $R'$ has fewer components
   than $R$, and its intersection matrix is still negative
   semi-definite, but not negative definite. It now
   follows by induction
   that there
   is a divisor as claimed, consisting entirely of components of
   $R'$.
\end{proof}

   We now give the

\begin{proof}[Proof of uniqueness]
   We claim first that in any
   decomposition $D=P+N$ satisfying the conditions of the
   theorem, the divisor $P$ is necessarily
   a \textit{maximal} nef
   $\Q$-subdivisor of~$D$. To see this, suppose that
   $P'$ is any nef
   divisor
   with $P\lessdivisor P'\lessdivisor D$.
   Then
   $P'=P+\sum_{i=1}^k q_iC_i$,
   where $C_1,\dots,C_k$ are the components of $N$
   and $q_1,\dots,q_k$ are
   rational numbers with
   $q_i\ge 0$.
   We have
   $$
      0\le P'\cdot C_j=\sum_{i=1}^k q_i\, C_i\cdot C_j
      \righttext{for } j=1,\dots k \ ,
   $$
   and hence
   $$
      \left(\sum_{i=1}^k q_iC_i\right)^{\!\!2}
      =\sum_{j=1}^k q_i \sum_{i=1}^k q_iC_i\cdot C_j\ge 0 \ .
   $$
   As the intersection matrix of $C_1,\dots,C_k$ is negative
   definite, we get
   $q_i=0$ for all $i$. So $P'=P$.

   To complete the proof it is now enough to show that
   a maximal effective nef $\Q$-subdivisor of $D$
   is in fact \textit{unique}.
   This in turn follows from:

   \begin{itemize}
   \item[(**)]
      If $P'=\sum_{i=1}^r x'_iC_i$ and
      $P''=\sum_{i=1}^r x''_iC_i$
      are effective nef $\Q$-subdivisors of $D$, then so is
      $P=\sum_{i=1}^r x_iC_i$, where $x_i=\max(x_i',x_i'')$.
   \end{itemize}
   As for (**):
   The divisor $P$ is of course an effective $\Q$-subdivisor of $D$,
   so it
   remains to
   show
   that it is nef, i.e., that the tuple
   $(x_1,\dots,x_r)$ satisfies the
   inequalities \eqnref{nefness-condition}.
   This, finally, is a consequence of the following elementary
   fact: Let $H\subset\Q^r$ be a halfspace, given by a linear
   inequality
   $
      \sum_{i=1}^r \alpha_i x_i\ge 0,
   $
   where the coefficients $\alpha_i$ are rational numbers with at
   most one of them negative. If two points $(x'_1,\dots,x'_r)$ and
   $(x''_1,\dots,x''_r)$ with $x'_i\ge 0$ and $x''_i\ge 0$ lie in $H$,
   then so does $(x_1,\dots,x_r)$, where $x_i=\max(x'_i,x''_i)$.
\end{proof}

\begin{varthm*}{Remark}\rm
   As experts may recognize, the maximality condition that
   defines $P$ is in the surface case equivalent to the defining
   condition of Nakayama's $\nu$-decomposition (see
   \cite[Sect.~III.1]{Nak04}). As Nakayama pointed out, it is
   also possible to obtain a proof by using results on
   $\nu$-decompositions and $\sigma$-decompositions (in
   particular \cite[Proposition~III.1.14]{Nak04},
   \cite[Lemma~III.3.1]{Nak04}, \cite[Lemma~III.3.3]{Nak04}, and
   \cite[Remark~III.3.12 and the subsequent Remark~(1)]{Nak04},
   when combined with
   arguments making use of \cite[Lemma~7.3]{Zar62}
   and \cite[Lemma~7.4]{Zar62}.

   When viewed from the point of view of $\nu$-decompositions,
   the essential content of the present note is to provide a
   quick, simple, and self-contained
   proof of the fact that in the surface case the
   $\nu$-decomposition of an effective $\Q$-divisor is a rational
   decomposition enjoying properties (ii) and (iii), and that it
   is the unique decomposition with these properties.
\end{varthm*}

\paragraph{Acknowledgements.}
   The author was partially supported by DFG grant
   BA\,1559/4-3. I~am grateful to F.~Catanese for helpful remarks
   and for pointing out a
   gap in the first version of this note.

\bigskip\footnotesize
   Tho\-mas Bau\-er,
   Fach\-be\-reich Ma\-the\-ma\-tik und In\-for\-ma\-tik,
   Philipps-Uni\-ver\-si\-t\"at Mar\-burg,
   Hans-Meer\-wein-Stra{\ss}e,
   D-35032~Mar\-burg, Germany.

\nopagebreak
   E-mail: \texttt{tbauer@mathematik.uni-marburg.de}

\end{document}